\documentclass[reqno,11pt,centertags, draft]{amsart}
\usepackage{amsmath,amsthm,amscd,amssymb,latexsym,upref}

\newcommand{\bbD}{{\mathbb{D}}}

\newcommand{\bbR}{{\mathbb{R}}}

\newcommand{\bbC}{{\mathbb{C}}}

\newcommand{\bbT}{{\mathbb{T}}}

\newcommand{\cK}{{\mathcal{K}}}

\newcommand{\cS}{{\mathcal{S}}}

\newcommand{\s}{\sigma}
\newcommand{\la}{\lambda}

\newcommand{\ov}{\overline}

\newtheorem{theorem}{Theorem}

\theoremstyle{definition}

\begin{document}

\title[Nevanlinna-Pick problem]
{On the Berg--Chen--Ismail theorem and the Nevanlinna-Pick
problem}

\author[L. Golinskii, F. Peherstorfer and P. Yuditskiy]
{L. Golinskii, F. Peherstorfer and P. Yuditskiy}

\address{Mathematics Division, Institute for Low Temperature Physics and
Engineering, 47 Lenin ave., Kharkov 61103, Ukraine}
\email{leonid.golinskii@gmail.com}

\address{Institute for Analysis, Johannes Kepler University Linz A-4040 Linz, Austria}
\email{franz.peherstorfer@jku.at, petro.yuditskiy@jku.at}
\thanks{
Partially supported by   the Austrian Founds FWF, project number:
P20413--N18 and {\it Marie Curie International Fellowship} within
the 6-th European Community Framework Programme, Contract
MIF1-CT-2005-006966}

\date{\today}

\keywords{moment problem, Blaschke product, Carleson measure, Pick
matrix} \subjclass{Primary: 30E05; Secondary: 30D50}

\begin{abstract}
In 2002 C. Berg, Y. Chen, and M. Ismail found a nice relation
between the determinancy of the Hamburger moment problem and
asymptotic behavior of the smallest eigenvalues of the corresponding
Hankel matrices. We investigate whether  an analog of this statement
holds for the Nevanlinna--Pick interpolation problem.
\end{abstract}

\maketitle

The Hamburger moment problem is a problem of finding conditions on a
sequence $\{s_j\}$, $j=0,1,\ldots$ so that there exists a positive
Borel measure $\s$ with infinite support and
\begin{equation}\label{hmp}
s_j=\int_{\bbR} x^j\,d\s(x), \qquad j=0,1,\ldots. \end{equation}
With $\s$ one can associate the infinite Hankel matrix and the
sequence of its principal submatrices
\begin{equation}\label{hm}
H=\|s_{i+j}\|_{i,j=0}^\infty, \qquad H_n=\|s_{i+j}\|_{i,j=0}^n,
\quad n=0,1,2,\ldots. \end{equation}
 By the famous result of Hamburger \eqref{hmp} has a solution if and
 only if $H_n\ge 0$ for all $n$.

 Denote by $\{\la_{n,j}\}_{j=0}^n$ the eigenvalues of $H_n$
\eqref{hm}. Because of the interlacing property we have
\begin{equation}\label{in}
0\le\la_{n+1,0}\le\la_{n,0}\le\la_{n+1,1}\le\la_{n,1}\le\ldots
 \le\la_{n+1,n}\le\la_{n,n}\le\la_{n+1,n+1}, \end{equation}
and so for each $k$ $\la_{n,k}$ are monotone decreasing,
$\la_k(H)=\lim_{n\to\infty}\la_{n,k}$ exist, and
\begin{equation}\label{lin}
0\le\la_0(H)\le\la_1(H)\le\ldots. \end{equation}

Let us call the problem \eqref{hmp} {\it regular} if $\la_0(H)>0$,
and {\it singular} otherwise. In 2002 C. Berg, Y. Chen and M. Ismail
\cite{bci} proved a beautiful result which states that \eqref{hmp}
is regular if and only if it has infinitely many solutions
(indeterminate).

The Hamburger moment problem is one of the representatives of the so
called classical interpolation problems \cite{akh}. In this note we
address to another one, specifically, the Nevanlinna--Pick
interpolation problem in the Schur class $\cS$ of functions
contractive and analytic in the unit disk $\bbD$.  This is a problem
of finding the solutions of
\begin{equation}\label{np}
f(z_k)=w_k, \qquad k=0,1,2,\ldots, \end{equation} where $z_k$ are
distinct points in $\bbD$, $w_k$ complex numbers, and $f\in\cS$. The
well known criterion for \eqref{np} to have at least one solution is
given in terms of Pick matrices by
$$ P_n:=\left\|\frac{1-w_i\bar w_j}{1-z_i\bar z_j}\right\|_{i,j=0}^n\ge 0$$
for all $n=0,1,\ldots$. For the eigenvalues $\{\la_{n,j}\}_{j=0}^n$
of the matrices $P_n$ the above relations \eqref{in}, \eqref{lin}
hold, and, again, we distinguish between the regular and singular
Nevanlinna--Pick problem, i.e., in this paper we say that the
problem \eqref{np} is regular if $\lambda_0(P)>0$, and it is
singular if  $\lambda_0(P)=0$.

With respect to a number of various questions there is a strong
similarity between different classical interpolation problems, that
is, if  one can prove this or that statement with respect to  one of
the classical problem a quite parallel statement holds for another
one. In this note we investigate the question: \textit{is it true
that a Nevanlinna-Pick problem has infinitely many solutions
(indeterminate) if and only if } $\lambda_0(P)>0$? We give a
negative answer to this question. More precisely, we construct data
$\{z_k,w_k\}_{k=0}^\infty$ of an indeterminate Nevanlinna-Pick
problem such that $\lambda_0(P)=0$.

First of all we note that the Blaschke condition on the
interpolation nodes $Z=\{z_k\}$
\begin{equation}\label{bla}
\sum_{k=0}^\infty (1-|z_k|)<\infty \end{equation} guarantees that
the interpolation problem
\begin{equation}\label{np0}
f(z_k)=0, \qquad k=0,1,2,\ldots, \end{equation} has infinitely many
solutions. Indeed, every function of the form
$$
f=g(z)B(z),
$$
where $g(z)\in \cS$, and $B(z)$ is the Blaschke product
$$
B(z)=\prod_k\frac{z_k-z}{1-z\bar z_k}\frac{|z_k|}{z_k}\,,
$$
solves \eqref{np0}. Note that \eqref{bla} is necessary and
sufficient for the problem \eqref{np0} to be indeterminate.

Thus, our goal is to construct a Blaschke set $Z$ such that
$\lambda_0(P)=0$ for the sequence of Nevanlinna-Pick matrices  of
the specific form $P_n=K_n$, where
$$ K_n:=\left\|\frac{1}{1-z_i\bar z_j}\right\|_{i,j=0}^n.$$
In fact, \textit{our main statement here characterizes completely
the regularity of such Nevanlinna--Pick problems  in the above
sense.}

Recall, see e.g. \cite{Garnett}, that a (finite) Borel measure $\nu$
on $\bbD$ is called a {\it Carleson measure}, if
\begin{equation}\label{cm}
\int_{\bbD} |f|^2\,d\nu\le C\int_{\bbT}|f|^2\,dm
\end{equation}
for all $f\in H^2$. Here $dm$ is the normalized Lebesgue measure on
the unit circle $\bbT$. Due to  Carleson's theorem such measures are
characterized completely by the following property: there exists
$C>0$ such that
$$
\nu(Q_\epsilon(\phi))\le C \epsilon
$$
for all $-\pi<\phi<\pi$, $0<\epsilon<1$,
$$
Q_\epsilon(\phi):=\{z\in\bbD: \phi-\pi\epsilon\le\arg
z\le\phi+\pi\epsilon, \ 1-\epsilon \le|z|\le 1\}.
$$

\begin{theorem}
Let $Z$ satisfy $\eqref{bla}$, $B(z)$ be the corresponding Blaschke
product. The Nevanlinna--Pick problem $\eqref{np0}$  is regular if
and only if the measure $\nu$, defined by
\begin{equation}\label{cm1}
 \nu(\{z_k\})=|B'(z_k)|^{-2},
\end{equation}
is a Carleson measure in $\bbD$.
\end{theorem}

\begin{proof}
Let $\nu$ \eqref{cm1} be the Carleson measure. For arbitrary
$c_0,\ldots,c_n\in\bbC$ put
$$  h(z)=B(z)\,\sum_{k=0}^n \frac{c_k}{z-z_k} \in H^2 $$
with
$$ \|h\|^2=\sum_{i,j=0}^n K_{ij} c_i\ov {c_j}, \qquad
K=\|K_{ij}\|_{i,j=0}^\infty=\left\|\frac1{1-z_i\bar
z_j}\right\|_{i,j=0}^\infty. $$ We have $h(z_j)=c_jB'(z_j)$ for
$j=0,1,\ldots,n$, and $h(z_j)=0$ for $j\ge n+1$, so
$$ \int_{\bbD} |h|^2\,d\nu=\sum_{j=0}^n |c_j|^2. $$
Hence by \eqref{cm}
$$ \int_{\bbT} |h|^2\,dm=\|h\|^2\ge \frac1{C}\,\sum_{j=0}^n |c_j|^2
$$
and so $\lambda_0(K)\ge C^{-1}>0$, as claimed.

Conversely, assume that the Nevanlinna--Pick problem in question is
regular. Put
$$\cK_B=(BH^2)^\bot=H^2\cap BH^2_-, \qquad
\phi_k(z)=\frac{B(z)}{z-z_k}, \quad k=0,1,\ldots. $$ It is easy to
see that $\{\phi_k\}$ is complete in $\cK_B$, so the system of
functions
$$
f(z)=B(z)\,\sum_{k=0}^n
\frac{c_k}{z-z_k}+B(z)g(z)=h(z)+B(z)g(z),\quad n=0,1,\ldots,
$$
$c_0,\ldots,c_n\in \bbC$, and $g\in H^2$, is dense in $H^2$. We
prove \eqref{cm} for such functions. As above,
$$
f(z_j)=c_jB'(z_j), \quad j=0,1,\ldots,n, \qquad f(z_j)=0, \quad j\ge
n+1 $$ and
$$
\int_{\bbD} |f|^2\,d\nu=\sum_{j=0}^n |c_j|^2, \qquad
\|f\|^2=\|h\|^2+\|g\|^2\ge \|h\|^2= \sum_{i,j=0}^n K_{ij} c_i\ov
{c_j}
$$
and so
$$
\int_{\bbT} |f|^2\,dm\ge \sum_{i,j=0}^n K_{ij} c_i\ov {c_j}\ge
\lambda_0(K)\,\sum_{j=0}^n |c_j|^2=\lambda_0(K)\,\int_{\bbD}
|f|^2\,d\nu. $$ The proof is complete.
\end{proof}

If the Nevanlinna--Pick problem \eqref{np0} is singular (as in the
example below), then so is the general problem \eqref{np}. Indeed,
for  $D={\rm diag}(w_0,w_1,\ldots,w_n)$
$$ P_n=K_n-DK_nD^*\le K_n, $$
and so $\la_0(P)\le\la_0(K)$.

{\bf Example}. Put $z_k=1-k^{-p}$, $p>1$. Evidently $Z$ is a
Blaschke set. On the other hand
\begin{align*}
    (1-|z_n|^2)|B'(z_n)| &=\prod_{k=1}^{n-1}\frac{z_n-z_k}{1-z_n z_k}
    \prod_{k=n+1}^{\infty}\frac{z_k-z_n}{1-z_n z_k} \\
    \le &\prod_{k=n+1}^{\infty}\frac{k^p-n^p}{n^p+k^p-1}
    \le \prod_{k=n+1}^{\infty}
    \left(1-\left(\frac n k\right)^p\right) \\ \le&
    \exp\left(-\sum_{k=n+1}^{\infty}\left(\frac n k\right)^p\right)\le
    \exp\left(-\frac{n+1}{2^{p}(p-1)}\right),
    \end{align*}
so
$$ |B'(z_n)|^{-2}\ge
\frac1{n^{2p}}\,\exp\left(\frac{n+1}{2^{p-1}(p-1)}\right). $$ Thus
the measure $\nu$ \eqref{cm1} is infinite and, moreover, it is not
of Carleson type.

\bigskip

There is a simple way of manufacturing regular Nevanlinna--Pick
problems \eqref{np0}. Recall that $Z=\{z_k\}$ is the Carleson
(uniformly separated) sequence if
\begin{equation}\label{car}
\delta(Z):=\inf_n \left|\prod_{k\not=
n}\frac{|z_k|}{z_k}\,\frac{z_n-z_k}{1-\bar z_n z_k}\right|>0.
\end{equation}
Assume that $Z$ satisfies \eqref{car}. By the theorem of H. Shapiro
and A. Shilds \cite{shs} the system of functions
$$ x_k(z)=\frac{(1-|z_k|^2)^{1/2}}{1-\bar z_k z} \quad k=0,1,\ldots
$$
forms the Riesz basis in $\cK_B$. So, for all $c_0, \ldots,
c_n\in\bbC$ there is $c>0$ such that
$$ \sum_{i,j=0}^n K_{ij}(1-|z_i|^2)^{1/2}\,(1-|z_j|^2)^{1/2}\,c_i\ov
{c_j}\ge c\,\sum_{j=0}^n |c_j|^2 $$ or
$$ \sum_{i,j=0}^n K_{ij}\,d_i\ov{d_j}\ge c\,\sum_{j=0}^n
\frac{|d_j|^2}{(1-|z_j|^2)}\ge c\,\sum_{j=0}^n |d_j|^2 $$ as
claimed.

\medskip

{\bf Acknowledgement}. The authors thank Alexander Kheifetz for
helpful discussions. The work is written mainly during  the first
author's visit to Johannes Kepler University, Linz. He wishes to
thank JKU  for the hospitality and the Marie Curie Foundation who
made this visit possible.

\end{document}